 \documentclass[letterpaper, 10 pt, conference]{ieeeconf}

\IEEEoverridecommandlockouts                              
\usepackage{placeins,graphicx,diagbox }

\usepackage{wrapfig} 
\usepackage{times} 
\usepackage{amsmath,tikz,pgfplots,mathtools,float}
\usepackage{amsmath}
\usepackage{subfigure}
\usepackage{multirow}
\usepackage{amssymb,latexsym}  
\usepackage{dsfont} 

\usepackage{amsmath,algorithm} 
\usepackage{amssymb}  
\usepackage{color}
\usepackage{cite}
\newcounter{algo}

\newtheorem{theorem}{Theorem}
\newtheorem{lemma}{Lemma}

\newtheorem{assumption}{Assumption}
\newtheorem{proposition}{Proposition}
\newtheorem{cor}{Corollary}
\newtheorem{remark}{Remark}

\def\uvs#1{{{\color{black}#1}}}

\newcommand{\wh}{\widehat}

\newcommand{\Real}{\mathbb{R}}

\newcommand{\Rscr}{\mathcal{R}}

\newcommand{\pmat}[1]{\begin{pmatrix} #1 \end{pmatrix}}

\newcommand{\epsilonbar}{{\bar \epsilon}}

 \newcommand{\remove}[1]{}

\def\Real{\mathbb{R}}

\def\argmin{\mathop{\rm argmin}}

\def\us#1{{{\color{black}#1}}}

\def\uvs#1{{{\color{black}#1}}}
\usepackage{verbatim}
\usepackage{ulem}
\title{Linearly Convergent Variable Sample-Size Schemes for Stochastic Nash Games:
 Best-Response Schemes and Distributed Gradient-Response Schemes}
\author{Jinlong Lei  and  Uday V. Shanbhag  \thanks{Lei and Shanbhag are with the
	Department of Industrial and Manufacturing Engineering, Pennsylvania
		State University, University Park, PA 16802, USA  {\tt\small
		jxl800,udaybag@psu.edu}. This research was partially
supported by the U.S. National Science Foundation grant CMMI-1538605 and
CAREER award CMMI-1246887 (Shanbhag).  A  extended version with proofs may be found on ArXiv and at \texttt{http://personal.psu.edu/vvs3/papers/CDC\_LS\_sub.pdf}}}

\date{}
\begin{document}

\maketitle
\thispagestyle{empty}
\pagestyle{empty}

\begin{abstract}
This paper considers an $N$-player stochastic  Nash game  in which  the $i$th   player
minimizes a composite objective $f_i(x) + r_i(x_i)$, where  $f_i$ is
expectation-valued and  $r_i$ has an efficient prox-evaluation. In this
context, we make the following contributions. (i)  Under a strong  monotonicity
assumption on the concatenated gradient map, we derive ({\bf optimal}) rate statements and oracle complexity
bounds for the proposed variable sample-size  proximal stochastic
gradient-response (VS-PGR) scheme; (ii) We overlay (VS-PGR) with a consensus
phase with a view towards developing distributed protocols for aggregative
stochastic Nash games. Notably, when the sample-size and the number of
consensus steps at each iteration grow at a suitable rate, a linear rate of
convergence can be achieved; (iii) Finally,  under a  suitable contractive
property  associated with the proximal best-response (BR)  map,   we design a
variable sample-size proximal BR (VS-PBR) scheme, where the proximal BR  is computed by
solving a sample-average problem. If the batch-size for computing the
sample-average is raised at a suitable rate, we show that the resulting
iterates converge at a linear rate and derive the oracle complexity.
\end{abstract}
\maketitle


\section{Introduction}
Noncooperative game theory~\cite{fudenberg91game,basar99dynamic} is a branch
of game theory \us{that considers} the resolution of conflicts among selfish players, each of
which tries to optimize \us{its} payoff function,   given the rivals'
strategies.  Nash games represent an important subclass of noncooperative
games, \us{originating from}  the  seminal work by~\cite{nash50equilibrium}. \us{Such models} have seen wide
applicability in a \us{breadth} of engineered systems, such as power grids,
communication networks,  and sensor networks.  In this paper, we consider   the
Nash equilibrium problem (NEP) with a finite set of  $N $ players indexed by
$i$  where $i \in\mathcal{N} \triangleq \{1, \cdots,N\}.$ For any $i \in
\cal{N}$, the $i$th player  \us{is characterized by} a strategy  $x_i \in \mathbb{R}^{n_i}$  and a
payoff function $ F_i(x_i,x_{-i})$  \us{dependent} on its strategy  $x_i$ and \us{parametrized
by} rivals' strategies  $x_{-i}\triangleq \{x_j\}_{j \neq i}   $. \us{If
$n \triangleq \sum_{i=1}^N n_i$} \us{and $x$ denotes the strategy profile, defined as} $x\triangleq (x_1,\cdots,x_N) \in
\mathbb{R}^n$.  \us{We consider  a}  stochastic \us{Nash game $\mathcal{P}$} \us{where the}  objective  of player $i$,  given rivals'  strategies
$x_{-i}$,  is to   solve the  following stochastic composite  optimization
problem:
\begin{align} \tag{$\mathcal{P}_i(x_{-i})$} \label{Ngame} \min_{x_i \in \mathbb{R}^{n_i}}  \
F_i(x_i,x_{-i})\triangleq f_i(x_i,x_{-i})+r_i(x_i) \end{align}
where $f_i(x
)\triangleq \mathbb{E}\left[\psi_i(x ;\xi(\omega)) \right],$   the random
variable  $\xi: { \Omega} \to \Real^d$  is defined on the probability space $({
\Omega}, {\cal F}, \mathbb{P})$,    $\psi_i : \mathbb{R}^n \times \mathbb{R}^d
\to \mathbb{R}$ is a scalar-valued  function, and
$\mathbb{E}[\cdot]$  denotes the expectation with respect to the probability
measure $ \mathbb{P} $.    We restrict our attention to nonsmooth convex Nash
games  where $f_i(x_i,x_{-i})$ is assumed to be smooth \us{and convex in $x_i$}
for any $x_{-i}$ \us{while} $r_i(x_i)$ is assumed to be convex \us{but a}
possibly nonsmooth function with an efficient prox-evaluation.  A  Nash
equilibrium (NE) of  the stochastic Nash game \us{in which the $i$th player solves the parametrized problem (${\mathcal P}_i(x_{-i})$)} is a \us{tuple}
$x^*=\{x_i^*\}_{i=1}^N \in \mathbb{R}^{n} $  such that  the following holds for
each player $ i  \in \mathcal{ N}$: \begin{align} F_i(x_i^*,x_{-i}^*) \leq
F_i(x_i ,x_{-i}^*)  \quad \forall x_i\in \mathbb{R}^{n_i} .  \nonumber
\end{align} In other  words,  $x^*$ is an NE if no player  can improve the payoff by
unilaterally deviating from the strategy $x_i^*$.

Our focus is two-fold: {\em (i) Development of variable
sample-size stochastic  proximal gradient-response (PGR) and proximal best-response (PBR) schemes with {\bf
optimal (deterministic)} rates of convergence; (ii) Extension of
PGR schemes to distributed (consensus-based) regimes, allowing
for  resolving aggregative games with a prescribed communication graph, where {\bf linear} rates of
convergence are achieved by combining increasing number of consensus steps with a
growing sample-sizes of sampled gradients.}

\noindent {\bf Prior research.} We discuss some relevant research on continuous-strategy Nash games and variance reduction schemes for stochastic optimization.

(i) {\em Deterministic Nash games.} Early work considered
convex Nash games (where players solve convex programs)
where the concatenated gradient map is either strongly
monotone~\cite{alpcan2002game} or merely monotone
maps~\cite{yin09nash2,kannan12distributed}.
While the aforementioned schemes utilized gradient-response techniques,
best-response schemes reliant on the contractive nature of the best-response
map were examined in~\cite{pang09nash}.

(ii) {\em Stochastic Nash games.} Regularized stochastic approximation schemes
were presented for monotone stochastic Nash
games~\cite{koshal13regularized} while extensions have been
developed to contend with misspecification~\cite{jiang18distributed} and the
lack of Lipschitzian properties~\cite{yousefian15selftuned}. More recently,
sampled best-response schemes have been developed in~\cite{pang2017two} while
rate statements and iteration complexity bounds have been provided for a class
of inexact stochastic best-response schemes
in~\cite{shanbhag16inexact,lei2017randomized,lei2017synchronous}. In fact, we
draw inspiration from our work in~\cite{lei2017synchronous} to develop superior
rate statements and extensions to distributed regimes.  Finally, a.s.  and
mean  convergence of sequences produced by BR schemes was proven in
\cite{lei2017randomized,lei2017asynchronous} for stochastic and misspecified
potential games.

(iii) {\em Consensus-based distributed schemes for Nash games.} Inspired by the
advances in consensus-based protocols for resolving distributed optimization
problems, Koshal et al.~\cite{koshal2016distributed} developed two sets of
distributed algorithms for monotone aggregative Nash games on graphs.  More
recently, in~\cite{yi2017distributed,belgioioso2017semi}, the authors combine
gradient-based schemes with consensus protocols to address generalized Nash games.

(iv) {\em Variance reduction schemes for stochastic optimization.} There has
been an effort to utilize increasing batch-sizes of sampled
gradients in stochastic gradient schemes, leading to improved rates of convergence,
as seen in strongly
convex~\cite{shanbhag15budget,jofre2017variance,jalilzadeh2018vss} and
convex regimes~\cite{jalilzadeh16egvssa,ghadimi2016accelerated,jofre2017variance,jalilzadeh2018vss}.

\noindent {\bf Novelty and Contributions.}

{(i). \em VS-PGR.}  In Section \ref{sec:monotone}, \us{under a strong
monotonicity assumption}, we \us{prove that} a variable sample-size  proximal
gradient response (VS-PGR) scheme is characterized by a linear rate of convergence in mean-suqared error
(Th.~\ref{thm1}) while in Th.~\ref{thm2}, we establish that  the iteration
complexity (in terms of proximal evaluations) and oracle complexity   to
achieve an $\epsilon-$NE denoted by $x$ where $x$  satisfies   $ \mathbb{E} [\| x-x^*\|^2 ] \leq
\epsilon$  are    $\mathcal{O} (\ln(1/\epsilon))$  and  $\mathcal{O} \left(
(1/\epsilon)^{1+\delta}\right)$, respectively, where $\delta\geq 0$ and
$\delta=0$ under a suitable selection of parameters.
Furthermore, it  is  shown in Corollary
\ref{cor1}  that with some specific   algorithmic parameters,
   the  iteration and oracle complexity    to obtain an $\epsilon-$NE are bounded by   $\mathcal{O} (\kappa^2\ln(1/\epsilon))$ and by  $\mathcal{O} \left( \kappa^2/\epsilon  \right)$,  where $\kappa$ denotes the condition number.

{(ii). \em Distributed VS-PGR.} In Section \ref{sec:agg},  addressing an open question in stochastic Nash game, we design  a
distributed  VS-PGR scheme  to compute an equilibrium of an aggregative  stochastic
Nash game over a communication graph. By increasing the number of consensus steps and sample-size at each iteration, this scheme is characterized by a linear rate of convergence  (Th.~
\ref{agg-thm1}). In Th.~\ref{agg-thm2}, we show that the iteration,
oracle, and communication complexity to compute an
$\epsilon$-Nash equilibrium are    $\mathcal{O} (\ln(1/\epsilon))$, $\mathcal{O} \left(
(1/\epsilon)^{1+\delta}\right)$, and $\mathcal{O} \left(\ln^2(1/\epsilon)\right)$ respectively, where  $\delta\geq 0$ and
$\epsilon-$NE$_2$ denotes an $x$ satisfying   $ \mathbb{E} [\| x-x^*\|  ] \leq \epsilon$.

 {(iii). \em VS-PBR.}  In Section \ref{sec:BR}, we develop a variable
sample-size proximal BR (VS-PBR) scheme (see Alg.~\ref {inexact-sbr-cont}) to
solve a class of stochastic Nash games with contractive proximal BR maps,
where   each  player solves a  \us{sample-average} best-response problem per
step. We  show in Th.~\ref{thm4}  that  the generated iterates converge to
the NE in mean at a linear rate under suitable  number of
scenarios, and also establish  that  the iteration  and oracle complexity   to
achieve an $\epsilon-$NE$_2$  are      $\mathcal{O}(
(\ln(1/\epsilon))$    and   $\mathcal{O} (
(1/\epsilon)^{2(1+\delta)})$ with  $\delta\geq 0.$

  {\bf Notation:} When referring to a vector $x$, it is assumed to
be a column vector while $x^T$ denotes its transpose. Generally, $\|x\|$ {denotes}
the Euclidean vector norm, i.e., $\|x\|=\sqrt{x^Tx}$.   We write \textit{a.s.} as the abbreviation for ``almost surely''.  For a real number $x $, we define by  $\lceil x \rceil $  the smallest integer greater than $x$.
 For a closed convex function    $r(\cdot)$, the proximal operator   is defined  in the  following for any   $\alpha>0$:
\begin{equation}\label{proximal}
\textrm{prox}_{\alpha r}(x) \us{\ \triangleq \ } \argmin_{y} \left( r(y)+{1\over 2\alpha} \|y-x\|^2\right) .
\end{equation}
 For simplicity,  $\xi$  denotes $\xi(\omega)$  and in a slight abuse of notation, $N$ denotes the number of players while $N_k$ denotes the batch-size of sampled gradients at iteration $k$.

\section{\us{Variable sample-size Gradient Response}}\label{sec:monotone}
This section considers \us{the development of a variable sample-size stochastic gradient response scheme for}  a class of    strongly  monotone   Nash games  associated with  a strongly  monotone  concatenated gradient map.  We
proceed to show that this scheme produces a sequence of iterates that converges
to the Nash equilibrium at a linear rate and establish the oracle complexity
to achieve an \us{$\epsilon$-Nash} equilibrium.

\subsection{Variable sample-size proximal \us{GR (VS-PGR)}}\label{sec:algorithm}
We impose the following assumptions   on $\mathcal{P}$.
\begin{assumption}~\label{assump-play-prob}  Let the following hold.\\
(a) The function $r_i$ is lower semicontinuous and convex with effective domain
denoted by $\Rscr_i \triangleq \mbox{dom}(r_i)$. Suppose $\Rscr \triangleq \prod_{j=1}^N \Rscr_i$ and $\Rscr_{-i} = \prod_{j \neq i} \Rscr_j$.\\
(b) $f_i(x_i,x_{-i})$ is \us{C$^1$} and convex in $x_i $ over on an open set
containing  $ \mathcal{R}_i$ for every fixed $x_{-i} \in \Rscr_{-i}$;\\
(c) For all  $x_{-i} \in \us{\Rscr_{-i}}$ and  any $\xi\in \Real^d$, $  \psi_i(x_i, x_{-i};\xi)$ is differentiable  in
$x_i$  over an open set containing    $\Rscr_i$  such that
 $ \nabla_{x_i} f_i(x_i,x_{-i})=  \mathbb{E}[\nabla_{x_i} \psi_i(x_i, x_{-i};\xi )]$.
\end{assumption}

If $ G(x;\xi) \triangleq \big( \nabla_{x_i}  \psi_i(x ;\xi ) \big)_{i=1}^N  $ and  $ G(x) \triangleq  \mathbb{E}[G(x ;\xi)],$
then   $ G(x ) = \left(    \nabla_{x_i} f_i(x)    \right) _{i=1}^N  $ by Assumption \ref{assump-play-prob}(iii).
The following lemma establishes a tuple $x^*$ is  \us{an NE of $\mathcal{P}$}  \us{if and only if it is a fixed point of a \us{suitable map}.} 
\begin{lemma}[{\bf Equivalence between NE and fixed point}] \label{lem1} Given   the stochastic Nash game $\mathcal{P}$, suppose
Assumption~\ref{assump-play-prob}  holds for each player $i\in \mathcal{N}$. Define $r(x) \triangleq ( r_i(x_i))_{i=1}^N.$
 Then $x^* \in X$ is an NE   if and only if  $x^*$   is a fixed point of  $\textrm{prox}_{\alpha r}(x -\alpha G(x ))$, i.e.,
\begin{align}\label{FP}
x^*=\textrm{prox}_{\alpha r}(x^*-\alpha G(x^*)), \quad  \forall \alpha>0.
\end{align}
\end{lemma}

%

\us{Suppose the iteration index is denoted by $k$  and player $i$'s strategy at
time $k$  is denoted by $x_{i,k} \in\mathbb{R}^{n_i}$}, which is an estimate of
its equilibrium strategy $x_i^*.$ We consider a \us{variable sample-size}
generalization of the standard  proximal stochastic gradient method, in which
$N_k$ sampled gradients are utilized at  iteration $k.$ \us{Given} a sample $\xi_{k}^1, \cdots ,\xi_{k}^{N_k}$ of $N_k$ realizations of the
random vector $\xi$,  for any $i\in \mathcal{N}$, given $x_{i,0}\in \us{\Rscr_i},$
player $i$ updates $x_{i,k+1}$ as follows:
$$x_{i,k+1}=\textrm{prox}_{\alpha  r_i} \Big[ x_{ i,k}-\alpha  { \sum_{p=1}^{N_k}  \nabla_{x_i}\psi_i(x_k;\xi_{k}^p)  \over N_k}\Big],$$
where  $\alpha>0 $ is the constant step size, $
\nabla_{x_i}\psi_i(x_k;\xi_{k}^p),~p=1,\cdots, N_k$ denote  the sampled
gradients. Define $w_k^p \triangleq
G(x_k;\xi_k^p)-G(x_k) $,  and $\bar{w}_{k,N_k}\triangleq{1\over N_k}
\sum_{p=1}^{N_k} w_k^p .$
Then the aforementioned scheme  can  be expressed as \begin{align} \tag{VS-PGR} \label{VSSA}
x_{k+1}=\textrm{prox}_{\alpha  r} \left[ x_k-\alpha \left( G(x_k )+\bar{w}_{k,N_k} \right) \right].
\end{align}
We make the following assumptions on the gradient map and the noise.
\begin{assumption}~\label{assp-noise}
(i) The mapping $G(x)$ is Lipschitz continuous over the set $  \Rscr $ with a constant $L$, namely,
$$\| G(x)-G(y)\| \leq L \|x-y\|\quad \forall x,y\in \Rscr.$$

\noindent
(ii) $G(x)$ is strongly monotone with parameter  $\eta,$
i.e.,
$$(G(x)-G(y))^T(x-y) \geq \eta  \|x-y\|^2\quad  \forall x,y\in \Rscr.$$

\noindent (iii) There exists a constant $\nu$ such that  the following holds for any $k\geq 0$:
    $\mathbb{E}[\bar{w}_{k,N_k} \mid \mathcal{F}_k]=0$,
 $\mathbb{E}[ \|\bar{w}_{k,N_k} \|^2\mid \mathcal{F}_k] \leq \nu^2/N_k\quad a.s.,$
where $\mathcal{F}_k\triangleq \sigma\{x_0,x_1,\cdots,x_{k}\}$.
\end{assumption}

\subsection{Rate analysis}

We begin with a simple  recursion  for  the conditional mean    squared error in  terms of sample size $N_k,$ step size $\alpha $,
and the problem parameters.

\begin{lemma}\label{lem2}
Consider   \eqref{VSSA} and let  Assumptions \ref{assump-play-prob} and \ref{assp-noise} hold.
Define $q\triangleq 1-2\alpha\eta +\alpha^2L^2.$ Then   for all $k\geq 0,$
\begin{align*}
\notag
\mathbb{E}[ \|x_{k+1}-x^*\|^2 |\mathcal{F}_k] &  \leq  q \|   x_k -x^* \|^2 + { \alpha ^2\nu^2 / N_k}, \quad a.s.
\end{align*}
\end{lemma}

Using Lemma \ref{lem2},  we   are able to  show the linear  convergence rate of algorithm \eqref{VSSA}.
\begin{theorem}[{\bf Linear convergence rate of VS-PGR}] \label{thm1}
Let     \eqref{VSSA} be applied to \us{$\mathcal{P}$}, where $N_k=\left \lceil \rho^{-(k+1)}  \right\rceil$ for some $\rho\in (0,1)$, and $\mathbb{E}[ \|x_0-x^*\|^2] \leq C$ for some constant $C>0.$
 Suppose Assumptions \ref{assump-play-prob} and \ref{assp-noise} hold, and  $\alpha  <2\eta/L^2 $. 
Define $q\triangleq 1-2\alpha\eta +\alpha^2L^2.$
Then  the following  holds for any $k\geq 0$.

\noindent (i) If $\rho\neq q,$ then $\mathbb{E}[ \|x_{k }-x^*\|^2]\leq
  C(\rho,q) \max \{\rho,q\}^k ,$ where $  C(\rho,q)  \triangleq  C+{\alpha^2 \nu^2 \over 1-\min\{\rho/q, q/\rho\}} .$

\noindent (ii) If $\rho= q,$ then for  any $\tilde{\rho}\in(\rho,1),$ $\mathbb{E}[ \|x_{k }-x^*\|^2]\leq \widetilde{D} \tilde{\rho}^k, $ where $\widetilde{D} \triangleq \left(C+{\alpha^2 \nu^2 \over \ln((\tilde{\rho}/\rho)^e)}\right)$.
\end{theorem}
%
%
\subsection{\us{Iteration and Oracle} Complexity} \us{Next, we examine the iteration (in terms
of proximal evaluations) and oracle complexity of this scheme to compute an
$\epsilon$-Nash equilibrium, defined next.} Recall that  a random strategy
profile $x: { \Omega} \to \Real^{n}$ is an $\epsilon-$NE if $\mathbb{E}[
\|x-x^*\|^2]\leq \epsilon$.
\begin{theorem}[{\bf Iteration and Oracle Complexity}] \label{thm2}
Let  \eqref{VSSA} be applied to $\mathcal{P}$, where $N_k=\left \lceil \rho^{-(k+1)}  \right\rceil$ for some $\rho\in (0,1)$, and $\mathbb{E}[ \|x_0-x^*\|^2] \leq C.$
 Suppose Assumptions \ref{assump-play-prob} and \ref{assp-noise} hold.  Define $q\triangleq 1-2\alpha\eta +\alpha^2L^2.$
Let $\alpha<2\eta/L^2$,  $\tilde{\rho}\in(\rho,1)$, $C(\rho,q)$ and $\widetilde{D}$ be defined in Theorem \ref{thm1}.
Then    the number of  proximal  \us{evaluations} needed to obtain an $\epsilon-$NE is bounded by $K(\epsilon)$, defined as
\begin{align}\label{rate-prox}
K(\epsilon) \triangleq \begin{cases}  &  {1\over \ln \left({1/q} \right)}  \ln \left({ C(\rho,q)\over \epsilon}\right)
\qquad \quad {\rm if ~} \rho<q<1,\\
&   {1 \over \ln \left(1/\tilde{\rho} \right)}
\ln\Big({ \widetilde{D} \over \epsilon}   \Big) \qquad  \qquad ~{\rm if ~} q=\rho,\\
& {1 \over \ln \left({1/\rho} \right)} \ln\left( { C(\rho,q)\over \epsilon}\right)  \qquad \quad {\rm if ~} q<\rho<1,
\end{cases}
\end{align}
and    the number of sampled gradients required is bounded by $M(\epsilon)$, defined as
\begin{align}\label{rate-SFO}
M(\epsilon) \triangleq \begin{cases}  &
 { 1  \over  \rho \ln (1/\rho)}  \left({ C(\rho,q)\over \epsilon}\right)^{\frac{\ln(1/\rho)}{\ln(1/q)}} +K (\epsilon) ~ {\rm if ~} \rho<q<1,\\&
{  1 \over  \rho \ln (1/\rho)} \left({\widetilde{D} \over  \epsilon} \right)^{\frac{\ln(1/\rho)}{\ln(1/ \tilde{\rho} )}} + K (\epsilon) \quad ~{\rm ~ if ~} q=\rho \\
& { 1  \over  \rho \ln (1/\rho)} \left({ C(\rho,q)\over \epsilon}\right)  +K(\epsilon)  \qquad
~ {\rm ~ if ~} q<\rho<1.
\end{cases}
\end{align}
\end{theorem}
{\bf Proof.} We \us{first} consider the case $\rho\neq q.$
  By Theorem \ref{thm1}(i),   the following holds:
$$\mathbb{E}[ \|x_{k }-x^*\|^2]\leq \epsilon ~\Rightarrow ~k\geq K_1(\epsilon) \triangleq {\ln\left({ C(\rho,q)/\epsilon}\right) \over \ln \left({1/ \max \{\rho,q\}} \right)}. $$
Then we achieve the    bound given in equation  \eqref{rate-prox} for cases $ \rho<q<1$ and $q<\rho<1$.
  Note   that for $\lambda > 1$ and any positive integer $K$, the following holds:
\begin{align}\label{exponent-sum}
\sum_{k=0}^K \lambda^k &  \leq \int_0^{K+1} \lambda^x dx \leq
\frac{\lambda^{K+1}} {\ln(\lambda)}.
\end{align}
 Then we may obtain  the following bound:
\begin{align*}    & \sum_{k=0}^{ K_1(\epsilon) -1}  N_k    \leq   \sum_{k=0}^{ K_1(\epsilon)-1}  \rho^{-(k+1)}   + K_1(\epsilon)
   \leq  {  \rho^{- K_1(\epsilon)}   \over  \rho \ln (1/\rho)} + K_1(\epsilon) .
\end{align*}
Note that for any $0<\epsilon, p<1, c_1>0$,  the following holds:
\begin{equation}\label{rel-exp}
\begin{split} & \rho^{-  \frac{\ln( c_1/\epsilon )}{\ln(1/p )} }
	  = \left(e^{\ln(\rho^{-1} )}\right)^{  \frac{\ln( c_1/\epsilon )}{\ln(1/p )} }
			 \\& =  e^{\ln(c_1/\epsilon ))^{\frac{\ln(1/\rho)}{\ln(1/p )}}}
			  	=   {(c_1/\epsilon )^{\frac{\ln(1/\rho)}{\ln(1/p )}}}  .
\end{split}
\end{equation}
Thus,  the number of sampled gradients required to obtain  an  $\epsilon-$NE  is   bounded by
$$  { 1  \over  \rho \ln (1/\rho)}
 \left( C(\rho,q)\over \epsilon\right)^{\frac{\ln(1/\rho)}{\ln(1/\max\{\rho,q\} )}} +K_1(\epsilon) .$$
Thus,  we achieve the    bound given in equation  \eqref{rate-SFO} for cases $ \rho<q<1$ and $q<\rho<1$.
The resultd for the case $\rho=q$  can be    similarly proved.
\hfill $\Box$

\begin{remark}\label{rem1}
(i) The above theorem establishes that the iteration  and oracle complexity to achieve an $\epsilon-$NE    are    $\mathcal{O} (\ln(1/\epsilon))$  and  $\mathcal{O} ( (1/\epsilon)^{1+\delta})$, where  $\delta=0$ when $\rho\in (p,1),$  $\delta= \frac{\ln(q/\rho)}{\ln(1/q)}~{\rm when ~} \rho<q<1$,
and $\delta=\frac{\ln(\tilde{\rho} /\rho)}{\ln(1/ \tilde{\rho} )}$ when $q=\rho.$
(ii) Suppose we use an alternative  metric to describe the $\epsilon$-Nash equilibrium: $x: { \Omega} \to \Real^{n}$ is an $\epsilon-$NE$_2$ if
$  \mathbb{E} [\| x-x^*\| ]  \leq \epsilon.$ Then  by  Jensen's inequality if follows  that
an $\epsilon^2$-NE is also an  $\epsilon$-NE$_2$, and hence by Theorem \ref{thm2} we obtain that
the iteration and oracle  complexity   to achieve an $\epsilon-$NE$_2$  are   $\mathcal{O} (\ln(1/\epsilon))$  and  $\mathcal{O} \left( (1/\epsilon)^{2(1+\delta)}\right)$, respectively.
\end{remark}

\begin{cor} \label{cor1}
Let   \eqref{VSSA} be applied to $\mathcal{P}$, where $\mathbb{E}[ \|x_0-x^*\|^2] \leq C$ for some constant $C>0.$
 Suppose Assumptions \ref{assump-play-prob} and \ref{assp-noise} hold.
Define the condition  number $\kappa\triangleq {L\over \eta}$.
Set $\alpha ={\eta\over L^2}  $
 and $N_k=\left \lceil \rho^{-(k+1)}  \right\rceil$  with $\rho=1-{1\over 2\kappa^2}$.
Then    the number of  proximal  \us{evaluations} and  samples required to obtain an $\epsilon-$NE \us{are} bounded by   $\mathcal{O} (\kappa^2\ln(1/\epsilon))$ and by  $\mathcal{O} \left( \kappa^2/\epsilon  \right)$, respectively.
\end{cor}
\section{Distributed VS-PGR for Aggregative Games}\label{sec:agg}
Next, \us{we consider an aggregative game $\mathcal{P}^{\rm agg}$, where the $i$th player solves the following parametrized problem:}
\begin{align} \tag{$\mathcal{P}^{\rm agg}_i(x_{-i})$} \label{Ngame_agg} \min_{x_i \in \mathbb{R}^{n_i}}  \
F^{\rm agg}_i(x_i,x_{-i}), \end{align}
where $F^{\rm agg}_i(x_i,x_{-i})\triangleq f_i(x_i,x_i+\bar{x}_{-i})+r_i(x_i)$,
 $\bar{x} \triangleq \sum_{i=1}^N x_i$ denotes the aggregate of all players'  decisions,  $\bar{x}_{-i}=\sum_{j=1,j\neq i}^N x_j$ denotes the aggregate of all players' decisions except player $i,$ and
 $f_i(x_i,x_i+\bar{x}_{-i}) \triangleq \mathbb{E}\left[\psi_i(x_i,x_i+\bar{x}_{-i} ;\xi ) \right] $    is expectation valued.  We impose the following assumptions   on the  stochastic aggregative game.
\begin{assumption}~\label{assump-agg1}  Let the following hold.\\
(a) The function $r_i$ is lower semicontinuous and convex with effective domain
denoted by $\mathcal{R}_i$  required to be  compact. \\
(b)   For any $y\in  \bar{\mathcal{R}}\triangleq \sum_{i=1}^N \mathcal{R}_i,$ $f_i(x_i,y)$ is C$^1$ and convex in $x_i $ over  an open set containing  $\mathcal{R}_i$.\\
(c) For all  $y\in\bar{\mathcal{R}}$ and  any $\xi\in \Real^d$, $  \psi_i(x_i, y;\xi)$ is differentiable  in
$x_i$  over an open set containing    $\mathcal{R}_i$  s.t.
 $ \nabla_{x_i} f_i(x_i,y)=  \mathbb{E}[\nabla_{x_i} \psi_i(x_i, y;\xi )]$.
\end{assumption}

\subsection{Algorithm Design}
We aim to design a distributed algorithm to \us{compute an NE of} $\mathcal{P}^{\rm agg}$, where each player may exchange information with its
local neighbors, and subsequently update \us{its estimate of the equilibrium strategy
and the aggregate.} The  \us{communication} among players  is \us{defined} by an
undirected  graph $\mathcal{G}=(\mathcal{N},\mathcal{E} )$, where  $\mathcal{N}
\triangleq\{1,\dots,N\}$  is the set of  players and  $\mathcal{E}$ is  the set
of undirected edges between players. The set of \us{neighbors of} player  $i$,
denoted $\mathcal{N}_i$, is defined as  $\mathcal{N}_i \triangleq \{ j\in \mathcal{N}:
(i,j)\in \mathcal{E}\}$.  Define  the adjacency matrix $ A=[a_{ij}]_{i,j=1}^N$,
where $a_{ij}>0$  if $j\in \mathcal{N}_i$ and $a_{ij}=0$, otherwise.    A
path in $\mathcal{G}$ with length $p$ from $v_1$ to $v_{p+1}$ is a  sequence of
distinct nodes, $v_1v_2\dots v_{p+1}$, such that  $(v_m, v_{m+1})\in
\mathcal{E}$, for all $m=1,\dots,p$. The graph $\mathcal{G}$ is termed {\it
connected} if there is a  path between  any two distinct players.

Though  each player does not have access to all players' decisions, it \us{may}
estimate   the aggregate  $\bar{x} $ \us{by communicating with its} neighbors.
Player $i$  at time $k$ holds an estimate $x_{i,k}$ for its equilibrium
strategy  and an estimate $v_{i,k}$ for the average of the aggregate.  To
overcome the fact that the communication network is sparse, we assume that to
compute $v_{i,k+1},  $    players communicate not once but $\tau_k$ rounds  at
major iteration $k+1 $. The   strategy of each player is updated by a variable
sample-size proximal stochastic gradient scheme   that depends on
parameters  $\alpha$  and $N_k,$ similar to   \eqref{VSSA} developed in Section
\ref{sec:monotone}. We now specify the scheme in Algorithm
\ref{algo-aggregative}.
\begin{algorithm}[htbp]
\small
\caption{Distrib. VS-PGR for Agg. Stoch. Nash Games} \label{algo-aggregative}
 {\it Initialize:} Set $k =0$, and  $ v_{\us{i,0}}=x_{\us{i,0}} \in \mathcal{R}_i$  for  any $i \in\mathcal{N}$.
Let $\alpha>0$ and $\{\tau_k,N_k\}$ be deterministic sequences.

{\it Iterate until convergence}

{\bf Consensus.} $\hat{v}_{i,k} := v_{i,k} ~  \forall i\in \mathcal{N}$ and
repeat $\tau_k$ times
$$\hat{v}_{i,k} := \sum_{j\in \mathcal{N}_i} a_{ij}\hat{v}_{j,k} \quad \forall i\in \mathcal{N}.$$
 {\bf Strategy Update.} for any $ i\in \mathcal{N}$

\begin{align}
x_{i,k+1} &:=\textrm{prox}_{\alpha  r_i}  \left[ x_{i,k}-\alpha \left(\nabla_{x_i}f_i(x_{i,k}, N\hat{v}_{i,k}) +e_{i,k}\right)\right], \label{alg-agg-strategy}
\\ v_{i,k+1} &:= v_{i,k}+x_{i,k+1}-x_{i,k},\label{alg-agg-average}
\end{align} where
$e_{i,k}    \triangleq {\sum_{p=1}^{N_k}
  \nabla_{x_i}\psi_i(x_{i,k}, N\hat{v}_{i,k};\xi_{k}^p) \over N_k}   -\nabla_{x_i}f_i(x_{i,k}, N\hat{v}_{i,k})  $.
\end{algorithm}

We impose the following conditions on the communication graph,   gradient mapping, and  observation noise.
\begin{assumption}~\label{assump-agg2}
(i) The undirected graph $\mathcal{G}$ is   connected  and the associated adjacency matrix $ A$ is \us{symmetric with row sums  equal to} one.

\noindent
(ii) $\phi(x)\triangleq  \left(    \nabla_{x_i} f_i(x_i,\bar{x})    \right) _{i=1}^N$ is strongly monotone over $\Rscr$ with parameter  $\eta,$
i.e.,
$(\phi(x)-\phi(y))^T(x-y) \geq \eta  \|x-y\|^2\quad  \forall x,y\in \mathcal{R}.$

\noindent (iii) The mapping $\phi(x)$ is Lipschitz continuous over $ \mathcal{R} $ with a constant $L$, i.e.,
$\| \phi(x)-\phi(y)\| \leq L \|x-y\|\quad \forall x,y\in \mathcal{R}.$

\noindent (iv)  For any $i\in \mathcal{N},$  $  \nabla_{x_i}f_i(x_i, y) $  is   Lipschitz continuous in $y$ over the set $\bar{\mathcal{R}} $ for every fixed $x_i\in \mathcal{R}_i$, i.e., there exists some constant  $L_i$ such that for  any $x_i\in \mathcal{R}_i$,
$$\|  \nabla_{x_i}f_i(x_i, y_1)- \nabla_{x_i}f_i(x_i, y_2)\| \leq L_i \|y_1-y_2\|\quad \forall y_1,y_2\in \mathcal{\bar{R}}.$$

\noindent (v) If $\mathcal{F}_k\triangleq \sigma\{x_0,x_1,\cdots,x_{k}\}$,
for any $i\in \mathcal{N},$ there exists a constant $\nu_i$ such that  the following holds for any $k\geq 0$:
    $$\mathbb{E}[e_{i,k}| \mathcal{F}_k]=0  {~~\rm and~~}
 \mathbb{E}[ \|e_{i,k}\|^2| \mathcal{F}_k] \leq {\nu_i^2 / N_k}\quad a.s.$$
\end{assumption}

\subsection{Convergence Analysis}
Define $ A(k)\triangleq  A^{\tau_k}$. Then by  Assumption  \ref{assump-agg2}(a),
$ A(k) $  is also  \us{symmetric} with the sum of each row \us{equaling} one. Note from the consensus step in
Algorithm~\ref{algo-aggregative} that $\hat{v}_{i,k}=\sum_{j=1}^N [ A(k)]_{ij} v_{j,k}$.  We  now recall a  prior result. 

\begin{lemma} \label{lem-pre}
  By Assumption \ref{assump-agg2} and \cite[Proposition 1]{nedic2009distributed}, there exists a constant $\theta>0$ and $\beta\in (0,1)$ such that
\begin{align}\label{geometric}
\left | \left[ A^k\right]_{ij}-{1\over N} \right| \leq \theta \beta^k\quad \forall i,j\in \mathcal{N}.
\end{align}

\end{lemma}

 We introduce the transition matrices $\Phi(k,s)$ from time $s$ to  time $k\geq s$ as follows:
  \begin{align*}
	\Phi(k,s) &= A(k) A(k-1)\cdots  A(s),\forall  ~ 0\leq s<k, \end{align*}
where $\Phi(k,k)  = A(k)$.
We  then obtain  the following  recursion on the mean-squared error.

\begin{proposition}\label{prp2}
Consider Algorithm \ref{algo-aggregative},
where $\tau_k=k+1$ and $N_k=\left\lceil  \beta^{-(k+1)/2} \right \rceil$.
Define   $M\triangleq \sum_{j=1}^N \max\limits_{x_j\in \mathcal{R}_j} \| x_{j}\|$,
$C_1\triangleq M\theta \Big(1+ 2  e\sqrt{1/\ln(\beta^{-1/2}) } \Big) $,   $C_2\triangleq { 4M \theta \over \ln(1/\beta)} ,$ and   $\varrho\triangleq \left(1-2\alpha\eta +2\alpha^2L^2\right) ,$ where  $\theta$ and $\beta$ are defined in Lemma \ref{lem-pre}.
Let Assumptions \ref{assump-agg1} and  \ref{assump-agg2} hold.  Then for any $ k\geq 0$,
\begin{equation}\label{agg-prp1}
\begin{split}
\mathbb{E}[\|x_{k+1}-x^*\|^2  ] & \leq \varrho \mathbb{E}[\|x_{k }-x^*\|^2  ]  +C_3\beta^{(k+1)/2},
\end{split}
\end{equation}
where $C_3$  is defined as \begin{align}\label{prp-defc}
C_3 & \triangleq  \alpha^2 \sum_{i=1}^N\nu_i^2+
4\alpha MN    \left(C_1  \beta^{1/2}+ C_2 \right) \sum_{i=1}^NL_i  \notag \\
	& +  4\alpha^2N^2 \left(C_1^2 \beta^{3/2}+ C_2 ^2 \beta^{1/2} \right)    \sum_{i=1}^NL_i^2.
\end{align}
\end{proposition}
{\bf Proof.} For purposes of brevity, we  merely outline the proof.
Firstly, we   give a recursion on  the conditional mean-squared error as follows: \vspace{-0.1in}
\begin{align*}
& \mathbb{E}[\|x_{k+1}-x^*\|^2 | \mathcal{F}_k] \leq \varrho \|x_{k}-x^*\|^2  +\alpha^2\sum_{i=1}^N\nu_i^2/N_k   \\
& +4\alpha MN \sum_{i=1}^NL_i   \|  \hat{v}_{i,k}-y_k\|   + 2\alpha^2N^2\sum_{i=1}^NL_i^2 \|  \hat{v}_{i,k}-y_k\|^2  .  \end{align*}\vspace{-0.1in}
We then establish an upper bound on the consensus error:
\begin{align*}
\|y_k-\hat{v}_{i,k}\|&\leq M\theta  \beta^{\sum_{p=0}^k \tau_p}  + 2  M \theta \sum_{s=1}^k  \beta^{\sum_{p=s}^k \tau_p}.\end{align*}
Finally, by getting an    upper bound  on  $\sum_{s=1}^k \beta^{\sum_{p=s}^k \tau_p}
\leq e  \sqrt{1/\ln(\beta^{-1/2}) } \beta^{ (k+1)(k+2)/2}  +{ 2\beta^{ (k+1) /2} \over  (k+1)\ln(1/\beta)} $ and taking the unconditional expectation,
we prove  the result. \hfill $\Box$

Based on Prop.~\ref{prp2},   we can   show  the  linear   rate of convergence of Algorithm \ref{algo-aggregative}.
\begin{theorem}[{\bf Linear convergence rate of Algorithm \ref{algo-aggregative}}] \label{agg-thm1}
Suppose  Assumptions \ref{assump-agg1} and  \ref{assump-agg2} hold. Consider Algorithm \ref{algo-aggregative},
where $\tau_k=k+1$, $N_k=\left\lceil  \beta^{-(k+1)/2} \right \rceil$ and   $\mathbb{E}[ \|x_0-x^*\|^2] \leq C$ for some  $C>0.$
 Let  $\alpha \in (0,\eta/L^2)$ and define  $\varrho\triangleq \left(1-2\alpha\eta +2\alpha^2L^2\right) .$
 Then we have the following assertions for any $k\geq 0$.

\noindent (i) If $\beta\neq \varrho^2,$ then $\mathbb{E}[ \|x_{k }-x^*\|^2]\leq
  \widetilde{C}(\varrho,\beta) \max \{\varrho,\sqrt{\beta}\}^k ,$ where $    \widetilde{C}(\varrho,\beta)   \triangleq  C+{C_3 \over 1-\min\{\varrho/\sqrt{\beta}, \sqrt{\beta}/\varrho\}} $ with $C_3$  defined in Proposition \ref{prp2}.

\noindent (ii) If $\beta=\varrho^2 ,$ then for  any $\tilde{\varrho}\in(\varrho,1),$ $\mathbb{E}[ \|x_{k }-x^*\|^2]\leq \widetilde{D}(\varrho) \tilde{\varrho}^k, $ where $\widetilde{D}(\varrho) \triangleq \left(C+{C_3 \over \ln((\tilde{\varrho}/\varrho)^e)}\right)$.
\end{theorem}

Similar to Theorem \ref{thm2}, \us{we may derive bounds on the iteration and
oracle complexity as well as the communication complexity to compute an $\epsilon$-Nash equilibrium.}
\begin{theorem} \label{agg-thm2}
Suppose  Assumptions \ref{assump-agg1} and  \ref{assump-agg2} hold. Consider Algorithm \ref{algo-aggregative},
where $\tau_k=k+1$, $N_k=\left\lceil  \beta^{-(k+1)/2} \right \rceil$ and   $\mathbb{E}[ \|x_0-x^*\|^2] \leq C$ for some constant $C>0.$
 Let  $\alpha \in (0,\eta/L^2)$ and define  $\varrho\triangleq \left(1-2\alpha\eta +2\alpha^2L^2\right) .$  Let  $\tilde{\varrho}\in(\varrho,1)$,
 $ \widetilde{C}(\varrho,\beta) $ and $\widetilde{D}(\varrho)$ be defined in Theorem \ref{agg-thm1}.
Then    the number of  proximal  evaluations  needed to obtain an $\epsilon-$NE is bounded as follows:
\begin{align*}
K(\epsilon) \triangleq \begin{cases}    {1\over \ln \left({1/\varrho} \right)}  \ln \left({ \widetilde{C}(\varrho,\beta) \over \epsilon}\right)
& {~\rm if ~} \beta<\varrho^2<1,\\
   {1 \over \ln \left(1/\tilde{\varrho} \right)}
\ln\Big({ \widetilde{D}(\varrho) \over \epsilon}   \Big) & ~{\rm if ~} \beta=\varrho^2,\\
 {1 \over \ln \left({1/\beta^{1/2}} \right)} \ln\left( { \widetilde{C}(\varrho,\beta )\over \epsilon}\right) &    {\rm if ~} \varrho<\beta^{1/2}<1,
\end{cases}
\end{align*}
and the round of communications is ${(K(\epsilon)+1)(K(\epsilon)+2) \over 2},$
and    the number of sampled gradients required is bounded by
\begin{align*}
M(\epsilon) \triangleq \begin{cases}
 { \left({  \widetilde{C}(\varrho,\beta)  \over \epsilon}\right)^{\frac{\ln(1/\beta^{1/2} )}{\ln(1/\varrho)}}  \over  \beta^{1/2} \ln((1/\beta^{1/2})}   +K (\epsilon)
 &  {\rm if ~} \beta<\varrho^2<1,\\
  { \left({\widetilde{D}(\varrho)
\over  \epsilon} \right)^{\frac{\ln(1/\varrho)}{\ln(1/ \tilde{\varrho} )}}   \over  \beta^{1/2} \ln((1/\beta^{1/2})}+ K (\epsilon) & {\rm if ~} \beta=\varrho^2, \\
   { \left({  \widetilde{C}(\varrho,\beta)  \over \epsilon}\right)  \over  \beta^{1/2} \ln((1/\beta^{1/2})}   +K(\epsilon)  &
{\rm if ~} \varrho<\beta^{1/2}<1.
\end{cases}
\end{align*}
\end{theorem}

\begin{remark}  Recall that in \cite{jakovetic2014fast}, a  fast distributed gradient algorithm  based on
 Nesterov's accelerated  gradient algorithm is employed  to solve a distributed convex optimization problem, where   at each  step,   $\mathcal{O}(\ln(k))$   consensus steps are taken. In~\cite{jakovetic2014fast}, the authors show that in merely convex settings,  the rate is $\mathcal{O}(1/k^2)$ (optimal)  and  total number of   communications rounds  is $\mathcal{O}(k\ln(k))$ up to time $k$.  Our scheme (Algorithm~\ref{algo-aggregative})  requires $\mathcal{O}(k^2)$ rounds of communications   to
recover the optimal linear rate of convergence but does so in a stochastic game-theoretic regime. In future work, we intend to investigate how the number of consensus   steps may be chosen to maintain geometric convergence while reducing  communication overhead.
\end{remark}

\section{Variable Sample-size Prox. Best Response}\label{sec:BR}
 In this section,  we consider the class of stochastic Nash games in which
the  proximal  BR map  is contractive~\cite{FPang09}.
We propose a  variable sample size   proximal BR scheme for computing an equilibrium, and derive  rate statements and establish iteration and oracle complexity bounds.
\subsection{Background  on proximal best-reponse maps}
For any $i\in \mathcal{N}$ and any  tuple $y\in \mathbb{R}^n,$  define the  proximal BR map $\wh{x}_i(y)$   as
$
\wh{x}_{i}(y)    \triangleq 
{\operatornamewithlimits{\mbox{argmin}}_{x_i \in \mathbb{R}^{n_i}}}
\left [  \mathbb{E}\left[ \psi_i(x_i,y_{-i};\xi ) \right]+ r_i(x_i) +{\mu \over 2} \|x_i-y_i\|^2 \right] $.
We impose the following assumption on problem \eqref{Ngame}.
\begin{assumption}~\label{assp-compact}
(i)   Assumption \ref{assump-agg1}(i).

\noindent
(ii)  For every fixed $x_{-i} \in \Rscr_{-i}$, $f_i(x_i,x_{-i})$ is \us{C$^2$} and convex in $x_i $ over on an open set
containing  $\mathcal{R}_i$.
\us{Moreover, $\nabla_{x_i} f_i(x_i,x_{-i})$ is assumed to be Lipschitz continuous  in $x_i$ uniformly in $x_{-i}$} with constant \us{$L_i$}, i.e., for any  $ x_i,x_i'\in\mathcal{R}_i,$
$$\|\nabla_{x_i}  f_i(x_i, x_{-i})-\nabla_{x_i}  f_i(x_i', x_{-i})\| \leq \us{L_i} \|x_i-x_i'\| .$$

\noindent (iii) For all  $x_{-i} \in \us{\Rscr_{-i}}$ and  any $\xi\in \Real^d$, $  \psi_i(x_i, x_{-i};\xi)$ is differentiable  in
$x_i$  over an open set containing    $\Rscr_i$.
Moreover,  for any $i \in \cal{N}$ and  all $x \in \mathcal{R} $, there exists   $M_i>0$ such that
  $\mathbb{E}[\| \nabla_{x_i}f_i(x)-\nabla_{x_i} \psi_i(x_i,x_{-i};\xi)\|^2] \leq  M_i^2.$
\end{assumption}

By Assumption \ref{assp-compact},  the \us{second derivatives}  of the   functions $f_i,~\forall i \in \mathcal{N}$ on $\mathcal{R}$ are bounded. Analogous to the avenue adopted in ~\cite{FPang09},   we may define  \begin{align}
&  \Gamma \triangleq \pmat{\frac{\mu}{\mu+\zeta_{1,\min}} &
	\frac{\zeta_{12,\max} }{\mu+\zeta_{1,\min}} & \hdots &
		\frac{\zeta_{1N,\max}}{\mu+\zeta_{1,\min}}\\
\frac{\zeta_{21,\max}}{\mu+\zeta_{2,\min}} &
	\frac{\mu}{\mu+\zeta_{2,\min}} & \hdots &
		\frac{\zeta_{2N,\max}}{\mu+\zeta_{2,\min}}\\
	\vdots & & \ddots & \\
\frac{\zeta_{N1,\max}}{\mu+\zeta_{N,\min}} &
	\frac{\zeta_{N2,\max}}{\mu+\zeta_{N,\min}} & \hdots &
		\frac{\mu}{\mu+\zeta_{N,\min}}} \label{matrix-hessian}
\end{align}
where $\zeta_{i,\min} \triangleq \inf_{x \in X} \lambda_{\min}  \left (\nabla^2_{x_i} f_i(x) \right)  \mbox{ ~and ~} \zeta_{ij,\max} \triangleq \sup_{x \in X}  \| \nabla^2_{x_ix_j} f_i(x) \|  ~  \forall j \neq i.$ Then by \cite[Theorem 4]{pang2017two} we obtain that
\begin{align}
  \pmat{\|\wh{x}_1(y') -\wh{x}_1(y)\| \\
				\vdots \\
\|\wh{x}_N(y') - \wh{x}_N(y)\|}  \leq   \Gamma   \pmat{\|y_1'- y_1\| \\
				\vdots \\
\|y'_N - y_N\|}  .\label{cont-prox-best-resp}
\end{align}
If the spectral radius    $ \rho(\Gamma)< 1$,  then  the  proximal best-response map is contractive w.r.t. some monotonic norm.
These sufficient conditions for the contractive property of the  BR   map $ \wh{x}(\bullet)$
can be found in~\cite{FPang09,pang2017two}.
\subsection{Variable sample-size proximal  BR scheme }
Suppose at iteration $k,$ we have  \us{$N_k$ samples} $\xi_{k}^1, \cdots
,\xi_{k}^{N_k}$ of the random vector $\xi.$ For any $x_i\in X_i,$ we
\us{approximate}  $f_i(x_i,y_{-i,k})$ by \us{${1\over N_k}\sum_{p=1}^{N_k} \psi_i(x_i,y_{-i,k}; \xi_{k}^p)$}
and  solve the  sample-average  best-response problem \eqref{SAA}.
We then get  the variable-size proximal  BR scheme (Algorithm~\ref{inexact-prox-br}).
\begin{algorithm}[!htbp]
\small
\caption{{Variable-size proximal  BR scheme}} \label{inexact-sbr-cont}
 Set $k:=0$. Let  $ y_{i,0}=x_{i,0} \in X_i$,   and  $\{\gamma_{i,k}\}_{k \geq 0}$
be a given  deterministic  sequence for $i = 1, \hdots, N$.
\begin{enumerate}
\item[(1)] For $i = 1, \hdots, N$,  player  $i$ updates estimate $x_{i,k+1}$ \us{as}
 \begin{equation}\label{SAA}
\begin{split}
x_{i,k+1}& =  \argmin _{x_i \in \us{\Real^{n_i}}} \Big[ {1\over N_k}\sum_{p=1}^{N_k} \psi_i(x_i,y_{-i,k}; \xi_{k}^p)
\\&\qquad \qquad+r_i(x_i)  + {\mu \over 2} \|x_i-y_{i,k}\|^2 \Big].
\end{split}
\end{equation}
\item[(2)] For $i = 1, \hdots, N$,  $y_{i,k+1} :=x_{i,k+1}$;
\item [(3)]  Set $k:=k+1$ and return to (1).
\end{enumerate}
\label{inexact-prox-br}
\end{algorithm}

\subsection{ Oracle and iteration complexity }

Define  $ \varepsilon_{i,k+1} \triangleq x_{i,k+1}- \wh{x}_{i}(y_k)$. We may obtain an bound on $\mathbb{E}[\|  \varepsilon_{i,k+1}\|^2]$ in the following lemma.
\begin{lemma}\label{lem4} Suppose Assumption   \ref{assp-compact}   holds.
Consider   Algorithm \ref{inexact-sbr-cont}.
Then  $ \mathbb{E}[\|  \varepsilon_{i,k+1}\|^2]\leq {M_i^2C_r^2 \over  N_k} ,$
where $C_r \triangleq {\mu \over \mu^2+L^2} (1-L/\sqrt{ \mu^2+L^2})^{-1}$ with  $L \triangleq \max_{i} L_i.$
\end{lemma}
{\bf Proof.}
Define $ \bar{w}_{i,k}(x_i)\triangleq {1\over N_k}\sum_{p=1}^{N_k}  \nabla_{x_i} \psi_i(x_i,y_{-i,k}; \xi_{k}^p)-\nabla_{x_i}f_i(x_i,y_{-i,k})$. By the optimality condition,  $x_{i,k+1}$ and  $\wh{x}_{i}(y_k) $ are respectively   a fixed point of the map ${\rm prox}_{\alpha r_i} \big[ x_i-\alpha\big( \nabla_{x_i}\tilde{f}_i(x_i,y_k)+ \bar{w}_{i,k}(x_i)\big)\big]$ and ${\rm prox}_{\alpha r_i} \big[ x_i-\alpha  \tilde{f}_i(x_i,y_k) \big] $ for any $\alpha>0$.  Then by the nonexpansive property of the
proximal operator,  the following  holds for any $\alpha>0:$
\begin{align*}
&\|\varepsilon_{i,k+1}\| \leq \alpha \|\bar{\uvs{w}}_{i,k}( x_{i,k+1})\|  \\&+ \Big\|  x_{i,k+1}- \wh{x}_{i}(y_k) -\alpha\left( \tilde{f}_i( x_{i,k+1},y_k)-  \tilde{f}_i(\wh{x}_{i}(y_k) ,y_k) \right)\Big\|  \\&\leq  \sqrt{ (1-\alpha\mu)^2+\alpha^2  L^2}  \|  \varepsilon_{i,k+1} \|+ \alpha  \|\bar{\uvs{w}}_{i,k}( x_{i,k+1})\|  .
\end{align*}
In the above inequality,  by setting  $\alpha={\mu \over \mu^2+L^2},$ we obtain that
$ \| \varepsilon_{i,k+1} \|\leq C_r \|\bar{\uvs{w}}_{i,k}( x_{i,k+1})\| .$
Then by using Assumption    \ref{assp-compact}(iii), the lemma is proved.
\hfill $\Box$

 Similar to    \cite[Prp. 4]{lei2017synchronous}, we can prove the linear rate of convergence.
 We then      establish the iteration and oracle complexity  to obtain an  $\epsilon-$NE$_2$,
which is random strategy profile $x: { \Omega} \to \Real^{n}$ satisfies $  \mathbb{E} [\| x-x^*\| ]  \leq \epsilon.$

\begin{theorem} \label{thm4}
Suppose  Ass.~\ref{assp-compact}  holds  and     $a  \triangleq \| \Gamma\|<1$.
Let   Algorithm \ref{inexact-sbr-cont}  be applied to the  stochastic Nash game \eqref{Ngame},
 where $\mathbb{E}[\|x_{i,0}-x_i^*\|] \leq C $ and  $N_k=\left\lceil {\max_i M_i^2 C_r^2\over  \eta^{2k}}\right\rceil$ for some $\eta\in (0,1).$ Define $c\triangleq  \max\{a,\eta\}$, let $ \tilde{\eta} \in (c,1)$,  and $D= 1/\ln((\tilde{\eta}/c)^e)$. Then
the number of the deterministic optimization problems solved and samples  required by player $i$ to obtain  an  $\epsilon-$NE$_2$  are
$\mathcal{O}(\ln(\sqrt{N}/\epsilon ))$ and $\mathcal{O}\Big(\big(  {\sqrt{N}/ \epsilon }  \big)^{   2   \ln(1/\eta) \over \ln(1/ \tilde{\eta})}  \Big)$, respectively.
 \end{theorem}

\section{Concluding Remarks}\label{sec:conclusion}
We consider a \us{class of stochastic Nash games} where each player-specific
\us{objective} is a sum of an expectation-valued smooth function and  a convex
nonsmooth function. \us{We consider three schemes: (i) Variable sample-size
proximal gradient response (VS-PGR) for strongly monotone stochastic Nash
games; (ii) Distributed VS-PGR for strongly monotone aggregative Nash games;
and (iii) VS proximal best-response (VS-PBR) for stochastic
Nash games with contractive best response maps. Under suitable assumptions, we
show that all schemes generate sequences that converge at the ({\bf optimal})
linear rate and derive bounds on the computational, oracle, and communication
complexity.}

\def\cprime{$'$}

\end{document}